\documentclass[12pt]{article}
\usepackage{amssymb,amsmath}
\usepackage{hyperref}
\usepackage{graphicx}
\usepackage{color}

\begin{document}

\title{\LARGE\bf A simple pseudo-Voigt/complex error function}

\author{
\normalsize\bf S. M. Abrarov\footnote{\scriptsize{Dept. Earth and Space Science and Engineering, York University, Toronto, Canada, M3J 1P3.}}\, and B. M. Quine$^{*}$\footnote{\scriptsize{Dept. Physics and Astronomy, York University, Toronto, Canada, M3J 1P3.}}}

\date{July 25, 2018}
\maketitle

\begin{abstract}
In this work we present a simple approximation for the Voigt/comp-lex error function based on fitting with set of the exponential functions of form ${\alpha _n}{\left| t \right|^n}{e^{ - {\beta _n}\left| t \right|}}$, where ${\alpha _n}$ and ${\beta _n}$ are the expansion coefficients. The computational test reveals that the largest absolute differences for the real and imaginary parts of the complex error function are $0.037$ and $0.036$, respectively.
\vspace{0.25cm}
\\
\noindent {\bf Keywords:} Voigt function, Faddeeva function, complex error function, complex probability function, rational approximation, spectral line broadening
\vspace{0.25cm}
\end{abstract}

\section{Methodology description}

Suppose that the function ${e^{ - {t^2}}}$ can be expanded by a set of the damping functions of kind ${\alpha _n}{\left| t \right|^n}{e^{ - {\beta _n}\left| t \right|}}$ with corresponding expansion coefficients ${\alpha _n}$ and ${\beta _n}$. Then, we can write
\begin{equation}\label{eq_1}
{e^{ - {t^2}}} = \sum\limits_{n = 0}^N {{\alpha _n}{{\left| t \right|}^n}{e^{ - {\beta _n}\left| t \right|}} + {\epsilon _N}\left( t \right)},
\end{equation}
where ${\epsilon _N}\left( t \right)$ is the error term associated with corresponding integer $N$. As a simplest case we consider only first two terms in series expansion \eqref{eq_1}. Consequently, from equation \eqref{eq_1} it follows that
\begin{equation}\label{eq_2}
{e^{ - {t^2}}} = {\alpha _0}{e^{ - {\beta _0}\left| t \right|}} + {\alpha _1}\left| t \right|{e^{ - {\beta _1}\left| t \right|}} + {\epsilon _1}\left( t \right).
\end{equation}

We found empirically that the following set of expansion coefficients
$$
{\alpha _0} = 1, \,\, {\beta _0} = 5.5, \,\, {\alpha _1} = 5.5 \,\, \text{and}\,\, {\beta _1} = 2.75
$$
provides a reasonable fitting. Thus, we can rewrite equation \eqref{eq_2} as
\begin{equation}\label{eq_3}
{e^{ - {t^2}}} = {e^{ - 5.5\left| t \right|}} + 5.5\left| t \right|{e^{ - 2.75\left| t \right|}} + {\epsilon _1}\left( t \right).
\end{equation}

\begin{figure}[ht]
\begin{center}
\includegraphics[width=34pc]{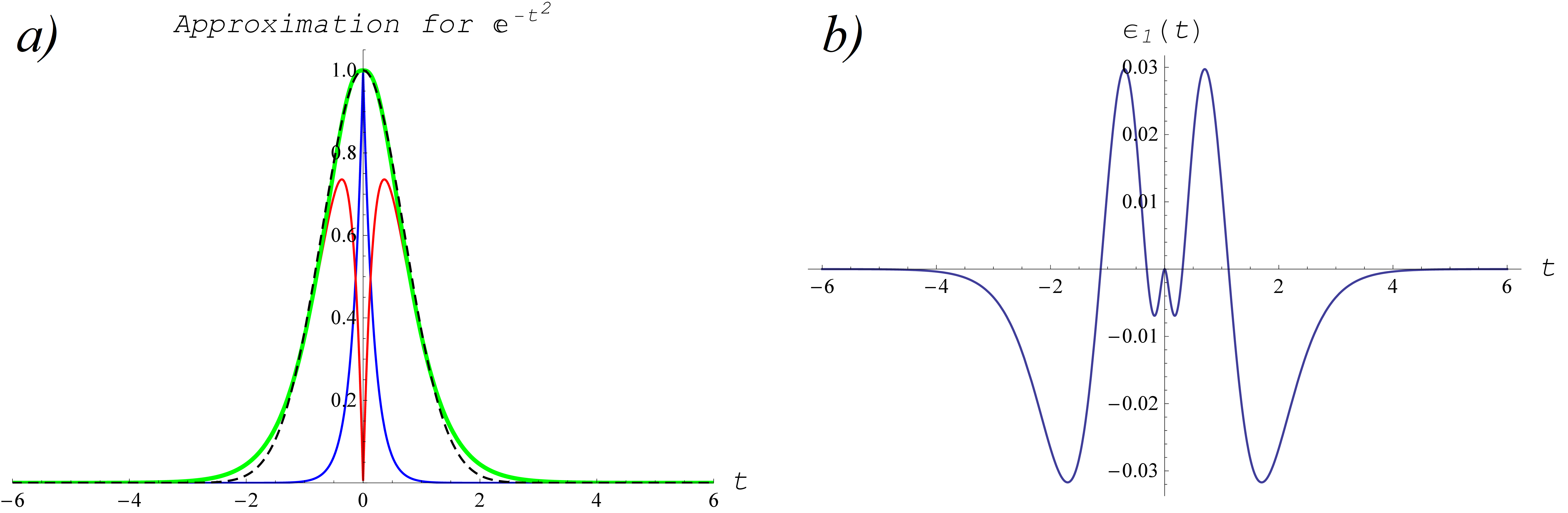}\hspace{2pc}%
\begin{minipage}[b]{28pc}
\vspace{0.5cm}
\small
{\sffamily {\bf{Fig. 1.}} a) The functions $e^{-5.5 \left|t\right|}$ (blue curve), $5.5\left|t\right|e^{-2.75\left|t\right|}$ (red curve), their sum $e^{-5.5 \left|t\right|} + 5.5 \left|t\right|e^{-2.75\left|t\right|}$ (green curve) and the original function $e^{-t^2}$ (dashed black curve). b) The error term $\epsilon_1 \left( t \right)$ defined as the difference ${e^{ - {t^2}}} - \left({e^{ - 5.5\left| t \right|}} + 5.5\left| t \right|{e^{ - 2.75\left| t \right|}} \right)$.}
\normalsize
\end{minipage}
\end{center}
\end{figure}

Figure 1a shows the functions ${e^{ - 5.5\left| t \right|}}$, $5.5\left| t \right|{e^{ - 2.75\left| t \right|}}$ and their sum ${e^{ - 5.5\left| t \right|}} + 5.5\left| t \right|{e^{ - 2.75\left| t \right|}}$ by blue, red and green curves, respectively. The original function  $e^{-t^2}$ is also shown by black dashed curve. Comparing the green and black dashed curves we can see that the function ${e^{ - 5.5\left| t \right|}} + 5.5\left| t \right|{e^{ - 2.75\left| t \right|}}$ approximates the original function ${e^{ - {t^2}}}$ reasonably well.

Figure 1b shows the error term ${\epsilon _1}\left( t \right)$. As we can see from this figure, the deviation is quite small since $\left|\epsilon _1\left( t \right) \right| << 1.$

\section{Derivation}

The complex error function, also known as the Faddeeva function, can be defined as \cite{Faddeyeva1961, Armstrong1967}
$$
w\left( z \right) = {e^{ - {z^2}}}\left( {1 + \frac{{2i}}{{\sqrt \pi  }}\int\limits_0^z {{e^{{t^2}}}dt} } \right)
$$
where $z = x + iy$ is the complex argument. It is not difficult to show that the complex error function \eqref{eq_3} can also be represented as given by
\begin{equation}\label{eq_4}
w\left( {x,y} \right) = \frac{1}{{\sqrt \pi  }}\int\limits_0^\infty  {\exp \left( { - {t^2}/4} \right)\exp \left( { - yt} \right)\exp \left( {ixt} \right)dt}.
\end{equation}
Consequently, separating the real and imaginary parts in equation \eqref{eq_4} results in
$$
w\left( {x,y} \right) = K\left( {x,y} \right) + iL\left( {x,y} \right)
$$
such that \cite{Armstrong1967, Srivastava1987}
\begin{equation}\label{eq_5}
K\left( {x,y} \right) = \frac{1}{{\sqrt \pi  }}\int\limits_0^\infty  {\exp \left( { - {t^2}/4} \right)\exp \left( { - yt} \right)\cos \left( {xt} \right)dt}
\end{equation}
and \cite{Srivastava1987}
\begin{equation}\label{eq_6}
L\left( {x,y} \right) = \frac{1}{{\sqrt \pi  }}\int\limits_0^\infty  {\exp \left( { - {t^2}/4} \right)\exp \left( { - yt} \right)\sin \left( {xt} \right)dt},
\end{equation}
respectively.

The real part $K\left( {x,y} \right)$ of the complex error function is known as the Voigt function \cite{Armstrong1967, Srivastava1987}. The imaginary part $L\left( {x,y} \right)$ of the complex error function has no specific name. Therefore, further we will refer to the function $L\left( {x,y} \right)$ as the $L$-function.

Define the constant $\gamma  = 2.75$. Then the equation \eqref{eq_3} can be expressed as
$$
{e^{ - {t^2}}} = {e^{ - 2\gamma \left| t \right|}} + 2\gamma \left| t \right|{e^{ - \gamma \left| t \right|}} + {\epsilon _1}\left( t \right).
$$
Making change of the variable $t \to t/2$ in this equation leads to
$$
{e^{ - {{\left( {\frac{t}{2}} \right)}^2}}} \approx {e^{ - \gamma \left| t \right|}} + \gamma \left| t \right|{e^{ - \frac{{\gamma \left| t \right|}}{2}}}
$$
or
\begin{equation}\label{eq_7}
{e^{ - {t^2}/4}} \approx {e^{ - \gamma t}} + \gamma t{e^{ - \frac{{\gamma t}}{2}}}, \qquad t \geqslant 0.
\end{equation}

Lastly, substituting the approximation \eqref{eq_7} into equations \eqref{eq_5} and \eqref{eq_6} yields
\begin{equation}\label{eq_8}
K\left( {x,y} \right) \approx \frac{1}{{\sqrt \pi  }}\left( {\frac{{y + \gamma }}{{{x^2} + {{\left( {y + \gamma } \right)}^2}}} + \frac{{4\gamma \left( {{{\left( {2y + \gamma } \right)}^2} - 4{x^2}} \right)}}{{{{\left( {4{x^2} + {{\left( {2y + \gamma } \right)}^2}} \right)}^2}}}} \right)
\end{equation}
and
\begin{equation}\label{eq_9}
L\left( {x,y} \right) \approx \frac{x}{{\sqrt \pi  }}\left( {\frac{1}{{{x^2} + {{\left( {y + \gamma } \right)}^2}}} + \frac{{16\gamma \left( {2y + \gamma } \right)}}{{{{\left( {4{x^2} + {{\left( {2y + \gamma } \right)}^2}} \right)}^2}}}} \right),
\end{equation}
respectively.

\section{Discrepancies}

In order to evaluate discrepancies it is convenient to define the absolute differences for the real and imaginary parts of the complex error function as
\[
{\Delta _{\operatorname{Re} }} = \left| {{K_{ref.}}\left( {x,y} \right) - K\left( {x,y} \right)} \right|
\]
and
\[
{\Delta _{\operatorname{Im} }} = \left| {{L_{ref.}}\left( {x,y} \right) - L\left( {x,y} \right)} \right|,
\]
respectively, where ${K_{ref.}}\left( {x,y} \right)$ and ${L_{ref.}}\left( {x,y} \right)$ are the references.

Figure 3a depicts the absolute difference for the pseudo-Voigt function \eqref{eq_8}. As we can see from this figure, the absolute difference ${\Delta _{\operatorname{Re} }}$ provided by the pseudo-Voigt function \eqref{eq_8} increases with decreasing parameter $y$. The largest discrepancy $0.037$ is observed at $y = 0.$ 

Figure 3b illustrates the absolute difference for the $L$-function approximation \eqref{eq_9}. Similar to the pseudo-Voigt function \eqref{eq_8}, the absolute difference of the $L$-function approximation \eqref{eq_9} increases with decreasing parameter $y$. The largest discrepancy $0.036$ also occurs at $y = 0.$

Although the proposed pseudo-Voigt function \eqref{eq_8} is not as accurate as that of reported in \cite{Ida2000}, it represents a very simple rational approximation (without any hyperbolic functions). Therefore, its application may be more convenient for rapid computation, for example, in debugging programs dealing with large-scale data.

\begin{figure}[ht]
\begin{center}
\includegraphics[width=34pc]{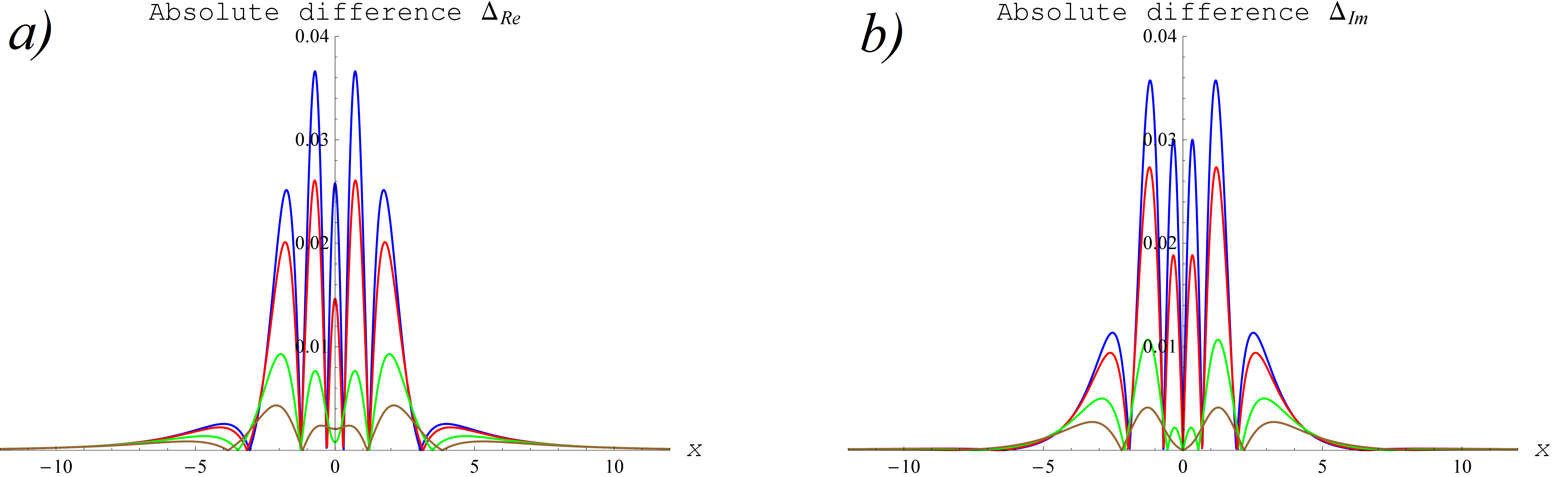}\hspace{2pc}%
\begin{minipage}[b]{28pc}
\vspace{0.5cm}
\small
{\sffamily {\bf{Fig. 2.}} Absolute differences for a) real and b) imaginary parts of the complex error function at $y = 0$ (blue curves), $y = 0.1$ (red curves), $y = 0.5$ (green curves) and $y = 1$ (brown curves).}
\normalsize
\end{minipage}
\end{center}
\end{figure}

\section{Conclusion}

A simple approximation for the Voigt/complex error function based on fitting with set of the exponential functions of form ${\alpha _n}{\left| t \right|^n}{e^{ - {\beta _n}\left| t \right|}}$ is presented. The largest absolute differences for the real and imaginary parts of the complex error function are found to be $0.037$ and $0.036$, respectively.

\section*{Acknowledgments}

This work is supported by National Research Council Canada, Thoth Technology Inc. and York University.

\bigskip

\end{document}